\documentclass[a4paper,12pt,leqno]{article}

\usepackage{geometry}
\usepackage{tikz-cd}\usepackage{bussproofs}
\usepackage{amssymb,amsmath}
\usepackage{proof}
\usepackage{psvectorian}
\usepackage{CJKutf8}
\usepackage[utf8]{inputenc}
\usepackage{graphicx}
\usepackage{graphicx}
\graphicspath{ {images/} }
\usepackage{qtree}
\usepackage{mathtools}
\usepackage{bussproofs}

\usepackage{hyperref}
\usepackage{listings}
\usepackage{indentfirst}
\usepackage{titlesec}
\numberwithin{equation}{section}

\newtheorem{defn}[equation]{Definition}

\newtheorem{rem}[equation]{Remark}
\newtheorem{exm}[equation]{Example}

\newtheorem{theorem}[equation]{Theorem}
\newtheorem{notat}[equation]{Notation}
\newtheorem{newpar}[equation]{}

\newtheorem{xdefn}{Definition.}
\newtheorem{xproposition}{Proposition.}
\newtheorem{xcorollary}{Corollary.}
\newtheorem{xrem}{Remark.}
\newtheorem{xexm}{Example.}
\newtheorem{xlemma}{Lemma.}
\newtheorem{xtheorem}{Theorem.}
\newtheorem{xnotat}{Notation.}
\newtheorem{xnewpar}{\it}
\newtheorem{xproof}{{\it Proof. }}
\newtheorem{xproofof}{{\it Proof}}

\newenvironment{definition}{\begin{defn}\em}{\end{defn}}

\newenvironment{example}{\begin{exm}\em}{\end{exm}}

\newenvironment{newparagraph*}[1]{\begin{xnewpar}\hspace*{-1.5mm}{#1}. \rm}{\end{xnewpar}}

\newenvironment{definition*}{\begin{xdefn}\em}{\end{xdefn}}
\newenvironment{remark*}{\begin{xrem}\em}{\end{xrem}}
\newenvironment{example*}{\begin{xexm}\em}{\end{xexm}}
\newenvironment{notation*}{\begin{xnotat}\em}{\end{xnotat}}
\newenvironment{proposition*}{\begin{xproposition}}{\end{xproposition}}
\newenvironment{corollary*}{\begin{xcorollary}}{\end{xcorollary}}
\newenvironment{lemma*}{\begin{xlemma}}{\end{xlemma}}
\newenvironment{theorem*}{\begin{xtheorem}}{\end{xtheorem}}

\titleformat*{\section}{\large\bfseries}

\begin{document}

\title{Stoic Logic and Natural Term Logic}
\author{Clarence Lewis Protin\footnote{Centro de Filosofia da Universidade de Lisboa.}}

\date{24th of November 2026}

\maketitle

\begin{abstract} In this paper we propose a reconstruction of Stoic Logic focused on the aspect of multiple generality using the Natural Term Logic (NTL) developped by the
author in \cite{ntl}. It can be shown, based on the available evidence, that Stoic logic was far more than a mere propositional logic (we cannot simply identiify the Stoic conditional with propositional conditionals). Its more general and complex logical theory - and notable the treatment of multiple generality - was rather patterned after the logico-syntactic mechanisms of natural language (subject to a regimentation which recalls Polish notation). NTL is likewise is a variable-free formal framework which captures the core syntactic and logic mechanisms of natural language (and notably a large array of features involving intensionality and anaphoric constructions). Bobzien and Shogry \cite{bobzienshogry} have argued that quantified conditionals and multiple generality were treated by the Stoics through a disambiguifying regimentation of syntax using indefinite pronouns and anaphoric constructions involving such pronouns.  We contend that the quantified condition expressed through an anaphoric construction involving indefinite pronouns  is the fundamental building block of ancient logic.  And we explore some possible building blocks for a formal reconstruction of the Stoic logic of universal conditionals (and its corresponding general form of Stoic modus ponens attested by ancient authors) based on NTL. More specifically, we present a version of NTL in which the anaphoric constructions of the Stoic quantified conditional correspond to a type of primitive constructor. The constructor also has a remarkable graph-theoretic interpretation which gives this version of NTL interest in its own right.  
\end{abstract}

\section{Introduction}

Both \cite{bobzienshogry} and \cite{sol} are dedicated to arguing that ancient logic was capable of expressing and reasoning about multiple generality. The first paper deals with Stoic Logic and the second focuses on an unorthodox reading of Aristotle as well as the testimonies of authors such as Galen and Boethius. The present work is dedicated to setting some foundations for interpreting the work \cite{bobzienshogry} in the light of a new version of Natural Term Logic (NTL)\cite{ntl}.

In Stoic logic the correspondence between the grammar (morphology and syntax) of natural language and the domain of meanings (\emph{lekta}) is not so straightforward. There are structural ambiguities and ellipses (cf. the so-called "donkey sentences") in natural language which renders it difficult to conclusively determine the associated \emph{lekton}.  However we must not fall into the opposite error (in interpreting ancient logic) of attempting an analysis of \emph{lekta} by discarding the structure of natural language altogether (as in modern predicate logic).  Rather, bearing in mind the reconstruction of a process of regimentation as outlined in \cite{bobzienshogry},  it seems desirable to employ a formal system which still matches closely the logical and syntactic mechanisms of natural language. This is where NTL comes in.

\section{Outline of Stoic Logic}

Sections 1 to 4 of \cite{bobzienshogry} provide a good introduction to Stoic logic accompanied by extensive documentary evidence for the reconstruction proposed. We also refer to this work for bibliographical references on Stoic logic. 

Stoic logic is a theory of what is meant or said (lekton) by natural language expressions (logoi).  It recognizes that the correspondence between the structure (grammar) of natural language expressions and the structure of sayables (lekta) is not perfect or unambiguous. This is why the Stoics sought to establish a regimentation of the syntax of natural language that would  render it viable as a mirror of the structure of the corresponding lekta.

The most basic distinction that can be made with regards to lekta is that between simple and compositie lekta and that between complete lekta (propositions, also called axiomata) and incomplete lekta. These last are further divided into kategoremata (monadic predicates) and less-than-kategoremata (dyadic or general polyadic predicates).  Composite lekta are defined by the Stoics inductively (just as in modern logic) in terms of simple lekta using certain constructors which include - but are no means limited to - what might be idenitied with modern propositiional connectives (classical or non-classical).  It appears -  at least in most cases - that  a single verb could correspond to a simple kategorema or less-than-kategorema depending on its valency.  And an example of a simple axioma (proposition) would be 'Socrates walks'. Apparently the Stoics also posited that the negation of such propositons counted as simple. Thus we can distinguish between affirmative and negative simple propositions.  However the authors in \cite{bobzienshogry} themselves admit - and adduce convincing evidence - that certain kinds of "conditional" involving anaphoric constructions (which we shall be mainly concerned with further ahead) were classified as simple.  This suggests that the inductive definition of lekta needs to be studied with greater care.  

Both Sextus Empricius and Diogenes Laertius give testimony that the Stoics divided affirmative simple propositions into definite (horismena), indefinite (aorista) and middle (mesa).  This corresponds to the subject to which is applied the predicate being a demonstrative pronoun, an indefinite pronoun or another kind of noun phrase (such a proper noun) respectively.  Thus indefinite pronouns were an essential component of lekta and, as we shall see, not only as subjects but also as part of the structure of kategorema and less-than-kategorema and formed the basis for the Stoic theory of quantification and multiple generality. The most basic type of definite affirmative simple proposition is "Socrates is walking". It can be noted that demonstrative propositions need not refer and that Galen in his work on the doctrines of Plato and Hippocrates gives testimony that Chryssipus considered the personal pronoun "I" as a demonstrative pronoun. 

The Stoic theory of inference (syllogism) was a sequent calculus (employing sequents with one consequent and two or more antecedents) defined by four rules (themata) and five (or possible seven) axioms (called anapodektoi or indemonstrables).  The Stoics distinguished between actual inferences and modes (tropoi) which are the analogue of the modern employment of variables to define the sequent-structure of arguments.  The indemonstrables express natural inferences relating to "implication" and exclusive disjunction and negation. The first one is usually interpreted as modus ponens $A, A\rightarrow B \vdash B$. As we shall see, there is reaason to doubt if that modern propositional formulation using a propositional connective does full justice to this Stoic inference.
According to the testimony of the late neoplatonic philosopher Simplicius, who lived in the 5th and 6th centuries CE, the third thema was nothing more than the cut rule (from $A,B \vdash C$ and $ C, D \vdash E$ we get $A,B, D \vdash E$) - where $D$ perhaps could be a multiset of formulae - while according to Apuleius the first thema was $A,B \vdash C \Rightarrow A, \sim C \vdash \sim B$. The third thema shows that Stoic sequents could have more than two antecedents. There is no direct evidence for the second and fourth themata though there have been many proposed reconstructions based on the indirect evidence of Sextus Empiricus, Alexander of Aphrodisias and Galen (in the context of the so-called Dialectical and Synthetic Theorems).

 According to Bobzien in her paper "Stoic Logic and Sequent Calculus"  the Stoic version of the sequence calculus was substructural and indeed a relevance logic.  However the reconstruction of the second and fourth themata by Dinucci and Duarte in their book "Introdução à Lógica Proposicional Estóica" (2016) suggests  rules which combine cut and contraction.  Dinucci and Duarte formulate the fourth thema as what amounts to: from $A,B \vdash C$ and $C, A, \Gamma \vdash D$ we get $A,B, \Gamma \vdash D$. Perhaps the fourth thema was in fact simply contraction $\Gamma, A, A \vdash B \Rightarrow \Gamma, A \vdash B$ and the second was a simpler case? Also there seems no doubt that the antecedents of Stoic sequents were multisets (i.e. the order of the formulae was immaterial).

\section{Natural Term Logic}

Natural Term Logic (NTL) is a formal system with the following aims:

\begin{enumerate}

\item To be formally similar to the grammatical structure of natural languagae. In particular

\begin{enumerate}
\item Not to make use of bound variables.
\item Be able to express propositional attitudes (intensionality) and "quantifying in".
\item Grammatically distinct expressions must correspond to distinct NTL terms (fine-grainedness).
\end{enumerate}

\item To come equipped with a system of formal transformation which transform a NTL term into a grammatically "simpler" term while preserving the core logical content.  We wish to have the situation in which two NTL term have the same core logical content if and only if they can be successively reduced to a common irreducibe term.

\end{enumerate}

The applications we have in mind for NTL are in the formalization of natural language, semantic networks, lexicology, philosophical logic, the philosophy of language, the formalization of Stoic logic,  type theory,  operad theory,  term-rewriting systems and complexity theory.

 Key examples of logical meaning invariance under grammatical transformations are the following:

  \begin{enumerate}
    \item (i) `Alice likes Bob' = `Bob is liked by Alice'.
  
    \item  (ii) `Alice likes Alice' =  `Alice likes herself'.
    
    \item  (iii) `Alice knows that Bob sleeps' = `Bob is known by Alice to sleep'.  

\item (iv)   'If something is a friend of Alice then it knows Bob` = `All Alice's friends know Bob'
   
  \end{enumerate}

Item (i) is an instance of the passive voice and  item (ii) an instance of the use of a  reflexive pronoun.  The left side of (iii) involves a propositional attitude.






 In the examples (i)-(iii) above not only do we intuitively assign the same basic logical meaning to the two sides but also the left-hand side is intuitively taken to be more simple or basic in some sense (the analogue of the `normal form' used in type theory): one could say that the left-hand sides render explicit the logical content implicit in the right-hand sides. NTL formalizes this intuition and setting up the corresponding reduction theory for NTL terms.

A \emph{weaver} $W$ of order $(n,m)$ is a surjective map $W: (1,2,...,n) \rightarrow (1,...,m)$ where $(1,2,..,k)$ denotes the ordered set of positive integers (which may be empty $()$).  A weaver of order $(n,m)$ can also be specified by an ordering of the elements of an equivalence relation on $(1,...,n)$. For example
the weaver $W$ of order $(3,2)$ given by $W(1) = 2$, $W(2) = 1$ and $W(3) = 1$ can be specified as $\{2,3\}\{1\}$. The name "weaver" comes from the fact that if we interpret the domain of a weaver as numbered and ordered strands then the weaver corresponds in general to an operation of permuting and joining these strands.

A \emph{selector} is a pair $(N,M)$ of sequences of non-negative integers of the same length such that the $i$th element of $M$ is smaller or equal than the $i$th element of $N$.

The language of NTL is defined as follows. We are given a collection of \emph{primitive} terms $A_1$, $A_2$,...,$A_n$,... including the special term $I$ and a finite collection of constructors $\Gamma_1,\Gamma_2,...,\Gamma_n$. Each primitive term $A$ has a \emph{valency} (saturation degree) which is a non-negative integer $s$. 
 Each constructor $\Gamma$ has a signature $(s_1,...,s_n)\rightarrow s_{n+1}$ where the $s_i$ are non-negative integers.

\begin{definition}
\emph{Terms} in NTL are defined as follows
\begin{itemize}
\item[-] A primitive term is a term.
\item[-] If $\Gamma$ is a constructor with signature $(s_1,...,s_n)\rightarrow s_{n+1}$ and $T_i$ are terms with valencies $s_i$ for $i=1,...,n$ then $\Gamma T_1...T_n$ is a term of valency $s_{n+1}$.
\item[-] Nothing else is a term.
\end{itemize}
\end{definition}

We sometimes indicate that a term $T$ has valency $s$  by writing $T^{(s)}$. The special term $I$ has valency $1$. \\

 We will also work only with the following set of constructors:

\[ \Pi^S, \Upsilon^W,  \Lambda^{(n,m)}_\wedge, \Lambda^n_\neg, \Lambda^n_\forall\]
where $S$ is a selector and $W$ is a weaver. We have that $\Pi^S$ has signature $(n, a_1,...,a_n) \rightarrow \Sigma b$ for $S = (a,b)$, $\Upsilon^W$ has signature $(n) \rightarrow m$ for $W$ of order $(n,m)$, $\Lambda^{(n,m)}_\wedge$ has signature $(n,m) \rightarrow n + m$ for $n,m\geq 0$,  $\Lambda^n_\neg$ has signature $(n) \rightarrow n$ and  $\Lambda^n_\forall$ signature $(n) \rightarrow n - 1$ for $n \geq 1$.
For a selector $S =(a,b)$ we will write $\Pi^SA B_1...B_n$  as $\Pi^bA B_1...B_n$ because $a$ is determined by the valencies of the $B_i$.  

  NTL expressions have the further important restriction that $I$ cannot occur as the argument of any constructor except $\Pi$ and in this case only in a non-head position with corrresponding entry 1 in the selector.


The valency of a NTL termt represents its ``degree of saturation", $0$ for  propositions and other kinds of objects ($0$-ary relations), $1$ for properties (or functions), $2$ for binary relations and so forth.  \\

The $\Upsilon$ family can be thought of as permuting the arguments of a term, as in the construction of the passive voice in $(i)$ of `Bob is liked by Alice', where a constructor $\Upsilon$ is applied to the term representing `likes' to yield `is liked by'. The $\Upsilon$ constructor represents also a generalized diagonalization of the arguments of a term, and in particular \emph{reflexivization} of a term: applying a certain $\Upsilon$ to `likes' yields the 1-valent term `likes oneself'. Thus if 2-valent primitive term $L$ represents "likes" then we have that the 2-valent term $\Upsilon^{\{2\}\{1\}}L$ represents "is liked by" and the 1-valent term  $\Upsilon^{\{1,2\}} L$ represents "likes oneself".\\

 The $\Pi$ family we be used to formalize both simple predication and more complex forms involving property formation seen in the discussion of (iii) . Let the 1-valent term $S$ represent "sleeps" and the 2-valent term $K$ represent "knows". let the 0-valent terms $A$ and $B$ represent "Alice" and "Bob". Then the most basic kind of predications are for instance:

\[ \Pi^{(0)} SA, \Pi^{(0,0)} LAB \]

These two 0-valent terms represent "Alice sleeps" and "Alice likes Bob".  Propositional 0-valent terms can also be arguments of these kinds of predications:

\[ \Pi^{(0,0)} KA \Pi^{(0)}SB\]

this 0-valent term represents the propositional attitude on the left side of example (iii) "Alice knows that Bob sleeps", where "Bob sleeps" is now an argument of "knows".

1-valent predicates can also be arguments of such types of predication. For example "to sleep is natural" which would be represented as $\Pi^{(0)}NS$ for the 1-valent predicate $N$ representing "being natural".  Also 1-valent predicates can modify 1-valent predicates (such as in adverbs and adjectives) to form new 1-predicates.  Let the 1-valent term $H$ represent "soundly".  Then the predicate "sleeps soundly"  is represented by the 1-valent term  $\Pi^{(1)} HS$ and the sentence "Bob sleeps soundly" by

\[ \Pi^{(0)} \Pi^{(1)} HS B \]

We can also represent instantiation of one (or more) of the arguments of binary (or $n$-ary) in general  predicates. For example from "knows" we form can properties "to be known by Alice" or "to know that Bob sleeps".  This is where the the special term $I$ comes in.  Here is the NTL representation of these last two  expressions:

\[ \Pi^{(0,1)} K A I,  \Pi^{(1,0)} K I \Pi^{(0)} S B\]

These are 1-valent terms.  We can also form "known by Alice to sleep" used in (iii):

\[ \Pi^{(0,1)} K A S \]

Thus "Bob is known by Alice to sleep" is :

\[ \Pi^{(0)} \Pi^{(0,1)} K A S B \]

Finally, the family $\Lambda$ includes logical connectives as well as a way of  expressing quantification without variables or variable binding. These logical constructors are completely general and thus can accommodate both classical and non-classical logics.

The examples (i)-(iii) above can be expressed in NTL as follows:

\begin{enumerate}
\item  $\Pi^{(0,0)}LAB = \Pi^{(0)} \Upsilon^{\{2\}\{1\}}LBA$
\item $\Pi^{(0,0)} LAA  = \Pi^{(0)} \Upsilon^{\{1,2\}}LA$
\item  $\Pi^{(0)} KA \Pi^{(0)}SB = \Pi^{(0)}\Pi^{(0,1)} KAS B$
\end{enumerate}

and the right side (iv) `all Alice's friends know Bob' can be expressed as\footnote{for simplicity we assume we have the constructor $\Lambda_\rightarrow$ which for the classical case could be defined in terms of the logical constructors for $\neg$ and $\wedge$.}
\[
\Lambda_\forall^1\Upsilon^{\{1,2\}} \Lambda^{1,1}_{\rightarrow}\Pi^{(0,1)}FAI\Pi^{(1,0)}KIB. \]

where $F$ represents "being  a friend".

Note the correspondence with the formulation the left side of (iv): the $\Upsilon$ constructor seems to mirror the anaphoric construction.

    We define a series of reductions on NTL terms which aim to capture meaning-preserving syntactic transformations ( transformations which preserved the basic logical meaning of a term) and each NTL term $T$ reduces to a unique normal term $N$.  The normal form of a term expresses - in a purely syntactic way - the core logical content of the term. But first we need some finite mathematics.

 Given two weavers $W$ of order $(n,m)$ and $U$ of order $(m,k)$ we obtain a weaver $UW$ of order $(n,k)$ by composition. Particular cases of weavers are weavers of order $(n,n)$ which are permutations - and  in particular the identity weaver given by $W(x) = x$.

For a weaver $W$ of order $(n,m)$ we denote by $W^\sharp$ the result of substituting each number $x$ in $(1,...,n)$ by $W(x)$. In the example above we have $W^\sharp = (2,1,1)$.  If $Z$ is a sequence of length $m$ then $W^\sharp Z$ is obtained by replacing each $y$ in $W^\sharp$ with $Z_y$, the $y$th element of $Z$. For example if $Z = ab$ then $W^\sharp Z = (b,a,a)$.

Consider a non-repeating sequence $S$ of length $m$ and a sequence $X$ of elements of $S$ of length $n$. Then we denote by $\mathcal{W}(X,S)$ the unique weaver of order $(n,m)$ such that  $(\mathcal{W}(X,S))^\sharp S = X$.

Given a two linearly ordered sets $A$ and $B$  of lengths $n$ and $m$ respectively and a surjective map $f:A \rightarrow B$ we can obtain a weaver $\mathcal{N}f : (1,...,n) \rightarrow (1,...,m)$ be defining $\mathcal{N}(x) $ as the position in $B$ of the image of the $x$th element of $A$.  

Given two weavers $W$ of order $(n,m)$ and $W'$ of order $(n',m')$ we define the sum $W+W'$  to be the weaver of order $(n + n', m + m')$ defined in the expected way by $W+W'(x) = W(x)$ for $x\leq n$ and $W+W'(x) = W'(x-n) + m$ otherwise. Given two linearly ordered sets $A$ and $B$ we denote by $A + B$ their concatenation starting with $A$ and ending with $B$.

Given the sequence $(1,...,n)$ an \emph{order partition} is a decomposition $(1,...,n) = I_1...I_k$ in which each $I_i$ is of the form $()$ or  $(a,...,b)$, the positive integers from $a$ to $b$ (including the case $(a)$) and in which all numbers of $I_{i+1}$ are larger than those in $I_i$. For example $(1,2)()(3)(4,5)$ is an order partition of $(1,...,5)$.  An order partition on $(1,...,n)$ can be specified equivalently by a sequence on non-negative integers $P= (p_1,...,p_k)$ such that $p_1 + ... + p_k = n$.  The last example could be specified by $(2,0,1,2)$. We also write $(1,...,n) / P = I_1...I_k$ for where $P$ is a sequence of non-negative integers of length $k$.  We will use the same notation when instead of $(1,...,n)$ we have any sequence of length $n$. For example $(abbcc) / (2,0,1,2) =  (ab)()(b)(cc)$.

Let $W$ be a weaver of order $(n,m)$ and $A$ be a subset of $(1,...,n)$ of cardinality $k$ with the induced order.  For example $(1,3,4)$ is an ordered subset of $(1,2,3,4,5)$. Then the restriction of $W$ to $A$, denoted by $W_A$ is
the weaver of order $(k,l)$ defined as follows. We consider the subset $W(A)$ of $(1,...,m)$ with the induced order and let $l$ be its cardinality. We define $W_A(x)$ follows. Let $y$ be the $x$th element of $A$ and let $z$ be the position of  $W(y)$ in $W(A)$.  Then $W_A(x) = z$.

Consider our example $W = \{2,3\}\{1\}$ and $A = (1,3)$. Then $W(A) = (1,3)$ and $W_A$ or order $(2,2)$ is given by $W_A(1) =2$ and $W_A(2) = 1$.

Given an order partition $P$ defined as $(1,...,n) = I_1...I_k$ and a weaver $W$ of order $(n,m)$ we say that $W$ is \emph{within} $P$  if $W = W_{I_1} + ... + W_{I_k}$. We say that  $W$ is \emph{without} $P$ if $W$ is not the identity weaver and its restriction $W_{I_i}$ to each $I_i$ of $P$ is the identity weaver. 
Let $W$ be a weaver or order $(n,m)$ within a partition $P = (I_1,...,I_k)$ on $A$. Then $W$ induces an order partition $P' = (J_1,...,J_k)$ on $(1,...,m)$.

An important result is that given an order partition $P$ defined as $(1,...,n) = I_1...I_k$ and a weaver $W$ of order $(n,m)$ $W$ can be factored in a unique way as  $W = W_{out}W_{in}$ where $W_{in}$ is a weaver within $P$ and $W_{out}$ is a weaver without the partition $P'$ induced by $W_{in}$.

Given a weaver $W$ of order $(n,m)$ and a sequence $S$ of non-negative integers of length $m$. Then a weaver $W \odot S$ is defined whose order depends on $W$. The construction is obtained by inserting in the place of $y$ in $(1,...,m)$ a fresh sequence $y^1,...,y^{S_y}$ where $S_y$ is the $y$th element of $S$, and for each $x$ in $W^{-1}(y)$ in $(1,...,m)$  we likewise insert in its place a fresh sequence   $x^1,...,x^{S_y}$. We then define a map $w$ between the resulting linearly ordered sequences $w: N \rightarrow M$ by $w(x^i) = (W(x))^i$.  We then take the corresponding weaver $\mathcal{N}w$ to be by definition $W \odot S$.

\begin{example} For example consider the weaver $L = ((1,2,3,4), (\{1,3\},\{2,4\})$ and $S = (2,2)$. Then we form first the map $f = ((1^1, 1^2, 2^1, 2^2, 3^1, 3^2, 4^1, 4^2), ( \{1^1, 3^1\}, \{1^2, 3^2\}, \{2^1, 4^1\} , \{2^2, 4^2\})$ and apply $\mathcal{N}$ to get
\[ L \odot S  = ((1, 2, 3, 4, 5, 6, 7, 8), ( \{1, 5\}, \{2, 6\}, \{3, 7\} , \{4,8\}) \]
\end{example}

 Given a sequence of integers $N$ we denote by $\Sigma N$ the sum of the elements of $N$. We consider that $\Sigma \epsilon = 0$.   Given a number $n$ and a non-negative ineger $M$ we denote by $n^M$ the sequence consisting of $i$ repetitions of $i$. If $i = 0$ this is considered to be the empty sequence.\\

 Given a selector $(N,S)$ with $|S| = k$ and a sequence of non-negative integers $T$  with $|T| = \Sigma S$.
Let 

\[T/S = \overline{\tau}_1 + ... + \overline{\tau}_k\]

Then we define the sequence $T/^N S$ by $(T /^N S)_i = \overline{\tau}_i + 1^{N_i - S_i}$ for $i = 1,...,k$ and

\[T/^\Sigma S =  (\Sigma \overline{\tau}_1,...,\Sigma\overline{\tau_k})\]

\begin{example}  Let $(N,S)$ be $((3,3),(2,1))$ and $T$ be $(2,2,1)$. Then  $(2,2,1) / (2,1) = (2,2) + (1)$. So 
$T /^N S = ((2,2,1),(1,1,1))$,  for we have $3-2 = 1$ and $3-1 = 2$, and $T /^\Sigma S = (4,1)$. 
\end{example}

We now define the reduction rules for NTL.

The central idea is to reduce the complexity of the head of a $\Pi$ application and to organize $\Upsilon$ applications in a canonical way. Reductions will allow us to define the concepts of normalization and normal form for NTL terms.  

Reductions will be divided into

\begin{itemize}
\item[-]  Structural reduction: the merging of sucessive $\Upsilon$-applications and dropping of trivial $\Upsilon$-applications.
\item[-] Predicative reductions, involving the removing of complex terms from the head of $\Pi$-applications and the simplification of trivial applications to series of $I$s.
\item[-] Pushing-in reductions.
\end{itemize}

The structural reductions consist in two rules:

\begin{equation}
\tag{$C_\Upsilon$} \Upsilon^{W_1} \Upsilon^{W_2} X \rightsquigarrow \Upsilon^{W_1W_2} X
\end{equation}

and

\begin{equation}
\tag{$Id_\Upsilon$} \Upsilon^{\{1\}...\{n\}} X^{(n)} \rightsquigarrow  X
\end{equation}

Predicative reductions have as goal to progressively simplify the head $T$ of a $\Pi$ application to the case in which the head is a primitive term. Thus we have reductions for each of the types of complex term $T$.

Let us start with the case in which $T$ begins with $\Pi$.

\begin{equation}
\tag{$R_{\Pi}$}
\Pi^A(\Pi^B TT_1...T_n)S_1...S_m \rightsquigarrow \Pi^{A /^\Sigma B} T U_1...U_n
\end{equation}

where $U_i = \Pi^{(A /^V B)_i}T_i\overline{S'}^i$ if $T_i$ is not $I$ and  $U_i = {S'}^i$ otherwise, where we use the selector $C = (V, B)$ determined by $B$ and the valencies of $T_1,...,T_n$, $S_1...S_m / B = \overline{S}^1 +...+\overline{S}^n$  and $\overline{S'}^i = \overline{S}^i I^{(V_i - B_i)} $ (the exponent means $V_i-B_i$ repetitions of the term $I$). If $(A /^V B)_i = \epsilon$ then $U_i = T_i$.

Assume that $T$ is a $\Upsilon^W$ application with $W$ a weaver. 

\[
\tag{$R_\Upsilon$}
\Pi^{A}(\Upsilon^W T)S_1...S_n \rightsquigarrow \Upsilon^{W\odot A} \Pi^{ W^\sharp A} T W^\sharp (S_1...S_n)\]

When the head $T$ is an $\Lambda^n_\neg$ application:

\[
\tag{$R_\neg$}
\Pi^{A}(\Lambda^n_\neg T)S_1...S_n  \rightsquigarrow \Lambda^{\Sigma A}_\neg  (\Pi^{A} T S_1...S_n).\]

The case in which $T$ is an $\Lambda^{n,m}_\wedge$ application is as follows. 

\[
\tag{$R_\land$}
\Pi^{A}(\Lambda^{n,m}_\land T^{(n)}S^{(m)})T_1...T_{n+m} \rightsquigarrow \Lambda^{\Sigma A_1, \Sigma A_2}_\land (\Pi^{A_1}TT_1...T_n)(\Pi^{A_2}ST_{n+1}...T_m)  \]

Here $A /(n,m) = A_1 + A_2$.


For $\Lambda^n_\forall$ applications:

\[
\tag{$R_\forall$}
\Pi^{A}(\Lambda^{n+1}_\forall T^{(n+1)})T_1...T_n \rightsquigarrow \Lambda^{1 + \Sigma A}_\forall (\Pi^{(1) + A} TIT_1...T_n)    \]

Then finally we have the predicative reduction: 

\[
\tag{$R_{I}$}
\Pi^{((1,1,...,1)} T II....I \rightsquigarrow T\]


Suppose we have a term of the form $\Upsilon^W\Pi^A TS_1...S_n$. Now $W$ can be decomposed as $W_{out}W_{in}$ for $W_{in}$ within the order partition $(1,2,...,\Sigma A) / A$. Consider the decomposition induced by the partition: $W_{in} = \Sigma_i W_i$ with $W_i$ of order $(o_i, p_i)$. Let $W$ be a weaver of ordern $(n,m)$. We use the notation $W+ 1^k$ to indiciate the sum of $W$ and the identity weaver of order $(k,k)$. Then we postulate the reduction

\[\tag{$P_\Pi$} \Upsilon^W\Pi^{A} TS_1...S_n \rightsquigarrow \Upsilon^{W_{out}} \Pi^{A'}T (\Upsilon^{W_1 +1^{V_1-A_1}}S_1)...(\Upsilon^{W_n + 1^{V_n-A_n}}S_n)             \]

where $A' =    (p_1,...,p_n)$ and $V_i$ is the valency of $S_i$.


Let $W$ be a weaver of order $(n,m)$ decomposed as $W_{out}W_{in}$ for the partition  $(1,2,...,n+m) / (n,m)$ with corresponding decomposition $W_{in} = W_1 + W_2$ of orders $(n,o)$ and $(m,p)$. Then we postulate:

\[ 
\tag{$P_\land$}
\Upsilon^W\Lambda^{n,m}_\wedge TS \rightsquigarrow \Upsilon^{W_{out}} \Lambda^{o,p}_\wedge (\Upsilon^{W_1}S)(\Upsilon^{W_2}T).   \]


The case of $\Upsilon^W \Lambda^n_\neg T$ for $W$ of order $(n,m)$ is:

\[ \tag{$P_\neg$} \Upsilon^W \Lambda^n_\neg T \rightsquigarrow \Lambda^m_\neg \Upsilon^W T . \]


Finally for $\Upsilon^W \Lambda^n_\forall T$ we postulate:

\[\tag{$P_\forall$} \Upsilon^W \Lambda^n_\forall T \rightsquigarrow \Lambda^{m+1}_\forall \Upsilon^{1+W} T  \]

where $W$ has order $(n-1,m)$. 

We note that the pushing-in reductions are in some sense redundant if the $\Upsilon^W$ application is a head of a $\Pi^A$ application. For the $\Upsilon^W$ can be eliminated by a predicative reduction. Pushing-in reductions are important  for instance for having a unique canonical version of an outermost $\Upsilon^W$ application of a term as well as a unique form for non-head arguments of $\Pi^A$- applications.  Pushing-in reductions are essential to obtain our result concerning reduction to a unique normal term.

We use the notation $T \rightsquigarrow^\ast S$ to indicate that $S$ is obtained from $T$ by multiple applications of reductions to subterms of $T$.

\begin{definition}
    An NTL term is called \emph{prenormal} if for all subterms of the form $\Pi T T_1...T_n$ we have that $T$ is a primitive term and there are no subterms of the form $\Pi^{ (1,1,...,1)} T I....I$ or $\Pi^{(1)} IT$.
\end{definition}

Note that applying only predicative reductions we always arrive in finitely many steps at a prenormal term. So prenormal terms do not allow further predicative reductions.


\begin{definition} A term of the form $\Upsilon^W T$ where $T$ starts with $\Pi$ or $\Lambda$ is called \emph{pushed-in} if the pushing-in reductions cannot be applied. 
\end{definition}

\begin{definition} A prenormal NTL term is called \emph{normal} if it all its subterms of the form   $\Upsilon^W T$ are pushed-in and the structural reductions cannot be applied. 
\end{definition}
\begin{theorem} Let $T$ be NTL term. Then there is a unique normal term $N$ such that $T\rightsquigarrow^\ast N$.

\end{theorem}

\begin{example} There is a tool that allows one to interactively reduce NTL terms until reaching a normal form.

For instance the NTL term 
\[\Pi^{(1,1,2)} \Upsilon^{\{1,4\}\{2\}\{3\}} \Pi^{(1,3)} \Pi^{(1,1)} C^{(2)} D^{(5)} D D D C B^{(1)} C\]

reduces to the normal term

\[\Upsilon^{\{1,2\}\{3\}\{4\}\{5\}} \Pi^{(1,4)} C \Pi^{(1,1,1,1,1)} D \Pi^{(1,1,1,1,1)} D C I I I I I I I I \Pi^{(4,1,1,1,1)} D \Upsilon^{\{4\}\{1\}\{2\}\{3\}\{5\}\{6\}} \Pi^{(1,2,1,1,1)} D B C C I I I I I I\]

The following file is an example of how the tool is used on the above term to reach the normal term above using several different reduction sequences.\\

\url{https://github.com/owl77/natlog/blob/main/ntltest3.txt}\\

\end{example}

Consider $\Pi^{(0)}\Pi^{(1)}L^{(1)} S^{(1)}J$.  This could be an adverbial construction: John is speaking loudly.  The 1-valent term "L" is $\Pi$-applied to the 1-valent term "S" to produce the 1-valent term "speaking loudly" which is then $\Pi$-applied to "John" to obtain a propositional 0-valent term. What would the reduction be? 
In this case we have $A = (0)$ and $B = (1,1)$ so that $A / B = (0)$ and so $A/^\Sigma B = (0)$ and $A/^V B = (0)$. 

\[ \Pi^{(0)}\Pi^{(1)}L^{(1)} S^{(1)}J  \rightsquigarrow          \Pi^{(0)} L \Pi^{(0)}SJ                       \] 

John is...speaking loudly  =  John is speaking...loudly.\\

Example (iii):

\[ \Pi^{(0)} \Pi^{(0,1)} K A S B \rightsquigarrow \Pi^{(0)} KA \Pi^{(0)} SB \]

\section{NTL and Stoic Conditionals}

As previously mentioned, Bobzien and Shogry \cite{bobzienshogry} have presented abundant evidence and arguments for the case that quantified conditionals and multiple generality were treated by the Stoics through a disambiguifying regimentation of syntax using indefinite pronouns and anaphoric constructions involving such pronouns.  We propose that the quantified conditional expressed through an anaphoric construction involving indefinite pronouns is the fundamental logical entity of ancient logic and that its closest modern equivalent is the $\Pi$ construction in dependent type theory (or Martin-Löf type theory). Our goal is to explore a reconstruction of the Stoic logic of universal conditionals based on NTL. Our strategy is to present a version of NTL in which the anaphoric constructions of the Stoic quantified conditional correspond to a type of primitive constructor. As we shall see, this constructor also has a remarkable graph-theoretic interpretation which gives this version of NTL interest in its own right.  

But first some general consideration on the basic interpretation of Stoic logic using NTL. We assume that we have a constructor $\Lambda_\rightarrow$  (for the propositional conditional) and  $\Lambda_\land$ which is to be read as Stoic exclusive disjunction. For simplicity we do not  distinguish between Chrysippian and Philonian conditionals. Both these new constructors have reductions corresponding analogous to $\Lambda_\wedge$.

The basis of our interpretation is that kategorema correspond to valency 1 SL terms. Axiomata correspond to valency 0 terms but not all valency 0 terms are propositional, some represent kind of objects which cannot be predicatively or functionally saturated (this will become clear later on). Polyadic kategoremata evidently correspond to SL terms of the corresponding valency (greater or equal than 2).  Primitive SL terms (of any valency) naturally correspond to simple lekta but, as remarked above,  we must be careful to distinguish with more care the nuances of the Stoic concept of "simple".
The basic form of predication which converts kategoremata or less-than-kategoremata into axiomata is captured by the constructor $\Pi^{(0,...,0)}AB_1...B_k$.  And also importantly the $\Pi$ constructor allows us to convert a less-than-kategorema into a kategorema. For intance if we have a term $D^{(2)}$  representing a less-than-kategorema than this can be converted into a kategorema via instationation through a term $S$ as follows:  $\Pi^{(0,1)} SI$ which has valency 1. Since SL uses Polish notation it is particularly apt to mirror the syntactical regimentation of Stoic logic which Bobzien and Shogry also associate to Polish notation.

The challenge is to formalize the regimented natural language syntax in which Stoic lekta were couched. In particular the use of indefinite pronouns, anaphora and the treatment of multiple generality.  Our approach
will be be focused on inference, specifically on versions of "modus ponens" and how they apply to various forms of Stoic "conditional sentences" which, as alreadys noted, do not alway correspond to modern propositional conditionals. Indeed, as we saw above, the authors in
\cite{bobzienshogry} argue there are kinds of conditionals which are not complex propositions and which cannot be "cut" (this technical term goes back to Chrysippus) into an antecedent and a consequent.  In its most basic form these read : if something is $A$ then it is $B$.  This is the simplest form of the Stoic "universal" conditional and also of the hypothetical syllogism (attested by Galen and Boethius) and which we argue is found in Aristotle as well.

In order to motivate the subsequent discussions, we now present a first naive approach to the universal conditional which involves in particular how the inference of the Stoic first indemonstrable (modus ponens) could be formalized in NTL and how reductions play a key role. Here for simplicity we omit the $judgement$ symbol $\vdash T$ expressing the assertion of 0-valent propositional terms $T$.

We postulate as rules modus ponens (the first indemonstrable):

\begin{prooftree}
\AxiomC{$\Lambda^{(0,0)}_\rightarrow A^{(0)}B^{(0)}$}
\AxiomC{$A$}
\BinaryInfC{$B$}
\end{prooftree}

and instantiation

\begin{prooftree}
\AxiomC{$\Lambda^1_\forall A^{(1)}$}
\AxiomC{$B^{(m)}$}
\BinaryInfC{$\Pi^{(0)} A B$}
\end{prooftree}

This rule correspond to $\forall$-elimination, and as it stands would not seem to have its direct counterpart in Stoic logic.
From these two rules we can derive the Stoic inferences from the universal conditional. An important example is the "conditional" in Aug. Dial.3.84 - 6, i. Fat. 11-15 (see \cite{bobzienshogry}): "If somebody is walking then he is moving". Then the Stoics had the inference from "Socrates is waking" to "Socrates is moving". The form of this inference is:

\begin{prooftree}
\AxiomC{$\Lambda^1\Upsilon^{\{1,2\}}\Lambda^{(1,1)}_\rightarrow W^{(1)}M^{(1)} $}
\AxiomC{$\Pi^{(1)}WS$}
\BinaryInfC{$\Pi^{(1)}MS$}
\end{prooftree}

Using NTL reductions and the two rules above it can be checked that this rule is admissable.

But now consider the important example in DL.7.75 discussed in \cite{bobzienshogry}: a clumbsy translation of the Greek would be "if something gives birth to something then that something is the mother of the other something". This sentence is then used as the premise in an inference of the form:

Based on these rules and logical equivalence under reduction of NTL terms we can derive the form of a

\begin{prooftree}
\AxiomC{$\Lambda^1_\forall \Lambda^2_\forall \Upsilon^{\{1,3\}\{2,4\}} \Lambda^{(2,2)}_\rightarrow A B$}
\AxiomC{$\Pi^{(0,0)}ACD$}
\BinaryInfC{$\Pi^{(0,0)}BCD$ }
\end{prooftree}

Instantiation of the first premise yields, using $D$:

\[ \Pi^{(0)} (\Lambda^2_\forall \Upsilon^{\{1,3\}\{2,4\}} \Lambda^{(2,2)}_\rightarrow A B) D \]
which reduces to
\[ \Lambda^{(1)} \Pi^{(1,0)} (\Upsilon^{\{1,3\}\{2,4\}} \Lambda^{(2,2)}_\rightarrow A B)) I D\]
 
and applying the rule again for $C$:
\[ \Pi^{(0)} (\Pi^{(1,0)} (\Upsilon^{\{1,3\}\{2,4\}} \Lambda^{(2,2)}_\rightarrow A B)) I D) C \]
which reduces to (exercise)

\[ \Lambda^{(0,0)}_\rightarrow (\Pi^{(0,0)} ACD)(\Pi^{(0,0)}BCD) \]

from which we get $\Pi^{(0,0)} BCD$ by modus ponens.

Our challenge is to investigate a primitive constructor NTL which is an alternative for the complex constructions $\Lambda^1\Upsilon^{\{1,2\}}\Lambda^{(1,1)}_\rightarrow W^{(1)}M^{(1)} $ and $\Lambda^1_\forall \Lambda^2_\forall \Upsilon^{\{1,3\}\{2,4\}} \Lambda^{(2,2)}_\rightarrow A B$ and which also expresses the fundamental structure of the use of analphora and relative pronouns in the corresponding natural language syntax. This is also motivated by the primitive $\Pi (x:A) B$ constructor in modern dependent type theory.

A key insight is that anaphoric constructions correspond to the "weaving" of the $\Upsilon$ constructor.

Consider $A$ and $B$ of valency 1. Then instead of $\Lambda_\forall \Upsilon^{\{1,2\}}\Lambda_\rightarrow AB$ we consider a single constructor $\Gamma^{\{1,2\}} AB$ in which - in the spirit of the geometric interpretation of \cite{ntl} - the strand of $A$ is  directly plugged into the input of  $B$.  We will discuss the graph calculus interpretation of $\Gamma$ further head. 
More generally given $A$ of valency $n$ and $B$ of valency $m$ we define the constructor $\Gamma^V AB$ where $V$ is a (possibly empty) set of ordered pairs $(a,b)$ in which $a \in \{1,...,n\}$ and $b \in \{1,...,m\}$ and such that for $(s_1,s_2),(t_1,t_2) \in V$ we have that $s_1 = t_1 \leftrightarrow s_2=t_2$. We call this set $V$ a $(n,m)$-linker. Such a linker defines a partial function $V: (1,...,n) \rightarrow (1,...,m)$ and we denote its domain (with the induced order) by $D(V)$ and the range by $R(V)$. The index of a $(n,m)$-linker is given by a pair $(p,q)$ where $p$ is $n - card(D(V))$ and $q$ is $m - card(R(V))$.  We define the function $\iota : (1,...,n) / D(V) \rightarrow \mathbb{N}$ such that $\iota(i)$ is the position of $i$ in $(1,...,n) / D(V)$. For example if $D(V) = \{(2,3)\}$ and $n = 4$ then $(1,...,4) / D(V) = (1,4)$ and $\iota(1) = 1$, $\iota(4) = 2$. We define $\kappa$ analogously but now for $(1,...,m)$ and $R(V)$.

We postulate that $\Gamma^V A B$ has valency $p + q$ where $(p,q)$ is the index of $V$. The linker is the direct interpretation of anaphoric constructions (cf. in certain grammatical theories anaphoric pronouns  are called "backwards pointing" pronouns and sometimes illustrated with similar lines). 
Thus we have a new version of NTL where $\Lambda_\forall$ has been replaced by $\Gamma$.  Future work will involve working out the reduction rules for $\Gamma$ and deriving the strong normalizability result of \cite{ntl}. We note that $\Gamma$ as a more satisfying graph theoretical interpretation compared to $\Lambda_\forall$ as presented in \cite{ntl}. Plain universal quantification could be recovered - as in dependent type theory - using a universe term $U^{(1)}$ for which we always have $\Pi^{(0)}UX$ for any term $X$. Thus $\Lambda_\forall X^{(1)}$ would be represented by $\Gamma^{(1,2)}UX$.
We now present the generalized modus ponens rule for the general case of $\Gamma$ and the particular instances relevant to Stoic logic.

\begin{prooftree}
\AxiomC{$\Pi^{(0,...,0)}\Gamma^V A^{(n)}B^{(m)}X_1....X_pY_{1}..Y_{q}$}
\AxiomC{$\Pi^{(0,...,0)} AZ_1....Z_n$}
\BinaryInfC{$\Pi^{(0,...,0)}B W_1....W_m$ }
\end{prooftree}

with the proviso that $Z_i = X_{\iota(i)}$ if $i \notin D(V)$ and that $W_i = Z_k$ if $V(k) = i$ and $W_i = Y_{\kappa(i)}$ otherwise.

This corresponds to the general case with parameters. An important specific instances is when the $(n,m)$-linker V has index $(n,m)$ which implies that $n = m$ and that all the threads are linked:

\begin{prooftree}
\AxiomC{$\Gamma^V A^{(n)}B^{(n)}$}
\AxiomC{$\Pi^{(0,...,0)} AZ_1....Z_n$}
\BinaryInfC{$\Pi^{(0,...,0)}B Z_1....Z_n$ }
\end{prooftree}

This instance includes the two kinds of generalized modus ponens we considered in the above discussion:

\begin{prooftree}
\AxiomC{$\Gamma^{(1,1)} A^{(1)}B^{(1)}$}
\AxiomC{$\Pi^{(0)} AZ$}
\BinaryInfC{$\Pi^{(0)} BZ$ }
\end{prooftree}

\begin{prooftree}
\AxiomC{$\Gamma^{(1,1)(2,2)} A^{(2)}B^{(2)}$}
\AxiomC{$\Pi^{(0,0)} AZ_1Z_2$}
\BinaryInfC{$\Pi B^{(0,0)} Z_1Z_2$ }
\end{prooftree}

We can also consider propositional modus ponens as an instance:

\begin{prooftree}
\AxiomC{$\Gamma^{\emptyset} A^{(0)}B^{(0)}$}
\AxiomC{$A$}
\BinaryInfC{$B$ }
\end{prooftree}

In Fig.1 we have a representation of the NTL terms $A$ and $B$ of valencies 4 and 3 respectively. Fig.2 represents the 3-valent term $\Gamma^V A B$ where $V = \{(3,2), (4,1) \}$. $V$ is a $(4,3)$-linker of index $(2,1)$. Fig 3 represents the term $\Gamma^{(1,1)}AB$ for $A$ and $B$ of valency 1. Figs.4 and 5 represent the last two generalized modus ponens rules above. Note how the rules have a geometric interpretation of the argument terms "sliding" (like a cable car) or being "transported" along the linking wires. Fig.6 represents a combination to two applications of the propositional modus ponens rule.

 \includegraphics[scale=0.4]{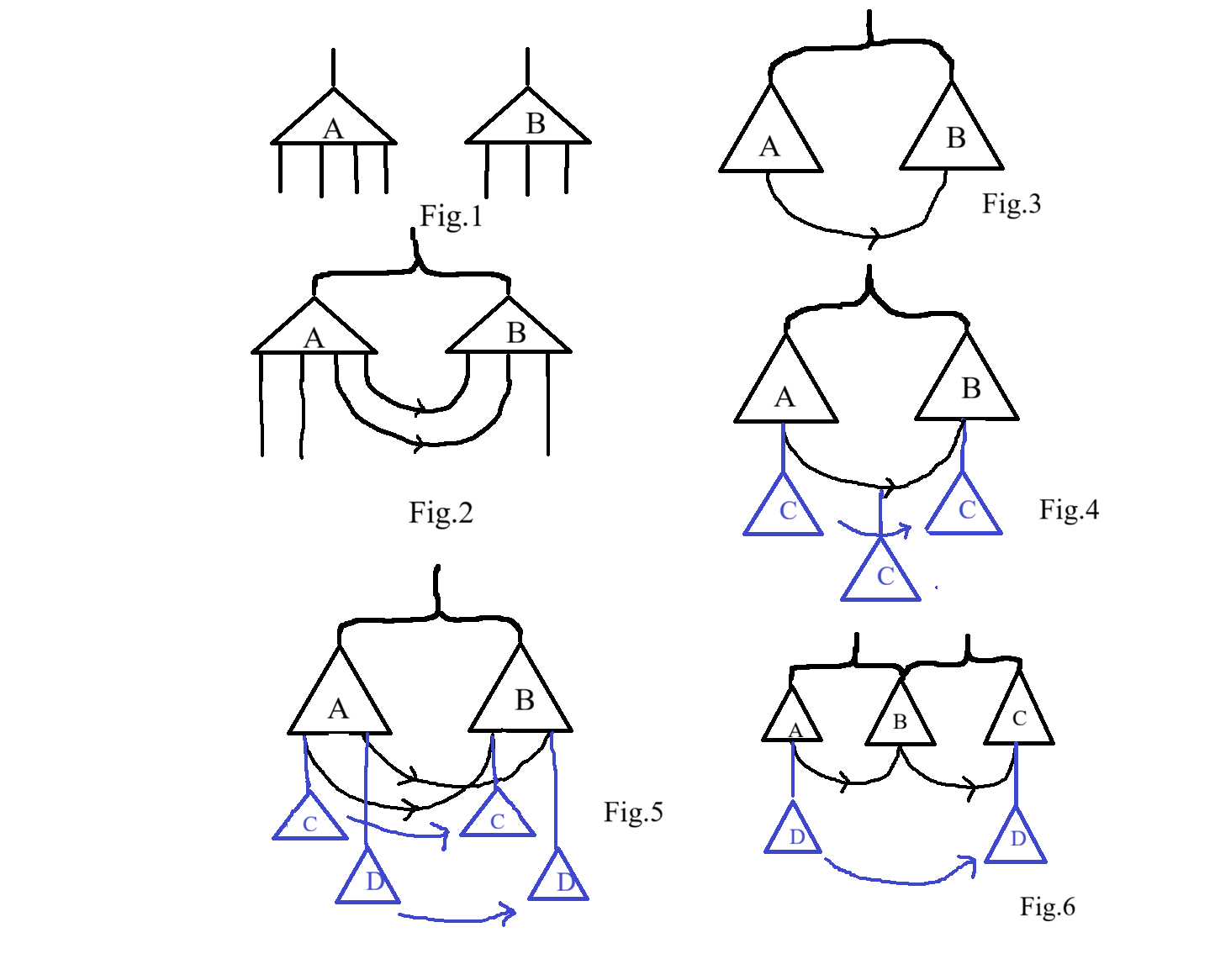}


\section{Conclusion}

We have explored some possible building blocks for a formal reconstruction of the Stoic logic of universal conditionals (and its corresponding general form of Stoic modus ponens attested by ancient authors) based a version of NTL in which the anaphoric constructions of the Stoic quantified conditional correspond to a type of primitive constructor. 
 The next step will be to explore a similar formalizations of the analogues of existential quantification (perhaps making some connections to the interpretation of ekthesis in \cite{sol} as well as $\Sigma$-types in dependent type theory). Then we can address the problem of reconstructing how reasoning about multiple generality may have been carrying out in Stoic Logic (which the authors in \cite{bobzienshogry} do not directly discuss.

Inspired by the treatment in \cite{sol} of natural deduction existential quantifier rules and the interpretation of \emph{ekthesis} in the context of Aristotle and Galen, we could inquire if Stoic logic was concerned with indefinite articles (with or without reference). Given a NTL term T or valency 1 we form the term $\epsilon T$ employing a new constructor $\epsilon$ of signature $(1) \rightarrow 0$ which represents "a T" which may or may not "exist".  Note that this is a case in which the resulting 0-valent term is not a propositional but an object.  In general this construction can be extended to terms of any valency and is related to the genitive construction(see \cite{sol} for details). If we were using the Peano operator corresponding to the definite article then in this case we would be in the presence of a function. In the present case we can think of the higher valency version of $\epsilon$ as the analogue of a "choice function".  Thus for $F^{(2)}$ we can interpret $\epsilon_2 F$ as a function which for a term $T$ chooses an element of $\Pi^{(0,1)}FTI$. That is to say $(\epsilon_2 F)T := \epsilon \Pi^{(0,1)}FTI$ which can be read "a $F$ of $T$".

We can combine $\epsilon$ with an existence predicate  (in the style of free logic) to obtain a constructor $\mathbb{X}$ applied to terms of any valency. For monadic predicates $T$ this reads "a T exists". For binary predicates $F$ the 1-valent term $(\mathbb{X}F)T$ reads " a $F$ of $T$ exists". Hence
"every man has a father" could be writen

\[ \Gamma^{(1,2)} M   (\mathbb{X} F) \]

We also should study possible interpretations of Stoic negation and perhaps take advantage of how NTL negation encompasses at once propositional and term negation.


\end{document}